\documentclass{article}

\usepackage{amssymb}
\usepackage{amsmath}
\usepackage{latexsym}

\newtheorem{definition}{Definition}[section]{\bf}{\it}
\newtheorem{theorem}[definition]{Theorem}{\bf}{\it}
{\bf}{\it}
\newtheorem{lemma}[definition]{Lemma}{\bf}{\it}
{\bf}{\it}
{\bf}{\it}
{\bf}{\it}

\def \In{\subseteq}
\def \IR {\mathbb R}
\def \IN {\mathbb N}
\def \IQ {\mathbb Q}

\def \s {\Sigma^*}
\def \om {\Sigma^\IN}
\def \dom {{\rm dom}}
\def \range {{\rm range}}

\def \ring {\mathcal R}
\def \algebra {\mathcal A}

\def \delm {\delta_{\mathcal M}}
\def \delmf {\delta_{{\mathcal M}_{<\infty}}}
\def \rmd {{\:\rm d}}

\def \pf {:\hspace{0.6ex}\subseteq \hspace{-0.4ex}}

\newcommand{\qq}{\phantom{.}\hfill$\Box$}
\newcommand{\pproof}{\noindent{\bf Proof:} }

\def \xoverline{\overline}

\title{\bf Computability of the Radon-Nikodym derivative}

\date{}

\author {Mathieu Hoyrup \thanks { \tt mathieu.hoyrup@loria.fr}\\
         LORIA - B248 615, rue du jardin botanique BP 239,\\
         54506 Vandoeuvre-l\`es-Nancy, France \\ \\
         Crist\'obal Rojas \thanks{\tt cristobal.rojas@unab.cl}\\
         Departamento de Matem\'aticas, Universidad Andres Bello, Chile\\ \\
         Klaus Weihrauch\thanks {\tt Klaus.Weihrauch@fernuni-hagen.de}\\
         Faculty of Mathematics and Computer Science, University of Hagen,\\
         Hagen, Germany
}

\begin{document}

\maketitle

\begin{abstract}
We study the computational content of the Radon-Nokodym theorem from measure theory in the framework of the representation approach to computable analysis. We define computable measurable spaces and canonical representations of the measures and the integrable functions on such spaces.
For functions $f,g$ on represented sets,  $f$ is W-reducible to $g$ if $f$ can be computed by applying the function $g$ at most once. Let $\rm RN$ be the Radon-Nikodym operator on the space under consideration and let $\rm EC$ be the non-computable operator mapping every enumeration of a set of natural numbers to its characteristic function. We prove that for every computable measurable space, $\rm RN$ is W-reducible to $\rm EC$, and we construct a computable measurable space for which
$\rm EC$ is W-reducible to $\rm RN$.
\end{abstract}


\section{Introduction}\label{seca}

The Radon-Nikodym theorem is one of the fundamental theorems in measure theory. In textbooks versions of various generality are proved
\cite{Bau72,Dud02,Rud87,Hal74,Bog07-1,Els99,Doo94}. We will study aspects of effectiveness of the following fairly general version.

\begin{theorem}[Radon-Nikodym]\label{t0}
Let $(\Omega,\algebra,\lambda)$ be a measured space where $\lambda$ is a $\sigma$-finite measure. Let $\mu$ be a finite measure that is absolutely continuous w.r.t. $\lambda$. Then there exists a unique function $h\in L^1(\lambda)$ such that for all $A\in\algebra$,
$$
\mu(A)=\int_A h\rmd \lambda.
$$
\end{theorem}

The function $h$ is called the {\em Radon-Nikodym derivative}, or {\em density} of $\mu$ w.r.t. $\lambda$, and is denoted by $\frac{\rm d\mu}{\rm d\lambda}$. \,
A measure $\mu$ is called {\em absolutely continuous} w.r.t. a measure $\lambda$, $\mu\ll\lambda$, if $(\lambda(A)=0\Longrightarrow \mu(A)=0)$ for all measurable sets~$A$.
 We mention that the condition
$(\forall A\in\algebra)\mu(A)=\int_A h\rmd \lambda$ is equivalent to the condition
$(\forall f\in L^1(\mu))\int f\rmd\mu=\int fh\rmd\lambda$.

In this article we ask whether the function $h$ can be computed from the measures $\lambda$ and~$\mu$. An answer can be given only relatively to computability concepts on the measures and functions under consideration, which must be defined in advance.
For studying computability on general spaces we use  the representation approach for computable analysis (TTE, Type Two theory of Effectivity)\cite{KW85,Wei00,BHW08}.
In TTE computability on finite words $w\in\s$ and infinite sequences $p\in\om$ is defined explicitly, for example by Turing machines, and then such finite or infinite sequences are used as names of abstract objects.

For the special application in measure theory we introduce computable measurable spaces and representations of measures and of integrable functions. In this setting, the Radon-Nikodym operator $\rm RN$
mapping $\lambda$ and $\mu$ to the function $h$ is not computable. We characterize its degree of non-computability in the $\leq_W$ hierarchy of problems on represented sets \cite{Wei92a,BG11,BG11a}. Let $\rm EC$ be the non-computable operator mapping every enumeration of a set of natural numbers to its characteristic function. We prove that for every computable measurable space the operator $\rm RN$ can be computed with a single application of $\rm EC$, that is, ${\rm RN}\leq _W{\rm EC}$. On the other hand we construct a simple computable measurable space such that $\rm EC$ can be computed with a single application of the operator $\rm RN$ on this space, that is, ${\rm EC}\leq_W{\rm RN}$.

In Section~\ref{secb} we summarize very shortly some definitions from TTE. In Section~\ref{sece}, we define computable measurable spaces with representations of the $\sigma$-finite measures and the finite measures and integrable functions.

 In Section~\ref{secc} we prove that for a computable measurable space the Radon-Nikodym operator $\rm RN$, mapping every $\sigma$-finite measure $\lambda$ and every finite measure $\mu\ll\lambda$ to the Radon-Nikodym derivative $h=\frac{\rmd\mu}{\rmd\lambda}$  can be computed with a single application of the operator $\rm EC$, that is ${\rm RN}\leq_W{\rm EC}$. In the proof we use two classical theorems: Levy's zero-one law \cite[Section~10.5]{Dud02} and the classical Radon-Nikodym Theorem.

In Section~\ref{secd} we construct a simple computable measurable space for which $\rm EC$ can be computed with a single application of the operator $\rm RN$, precisely ${\rm EC}\leq_{sW}{\rm RN}$
(see Section~\ref{secb}). Therefore, for sufficiently rich computable measurable spaces, ${\rm RN}$ and ${\rm EC}$ have the same degree of non-computability,  that is,    ${\rm RN}\equiv_W{\rm EC}$.

In Section~\ref{secf} we discuss other proofs of  ${\rm RN}\leq_W{\rm EC}$ and give a condition sufficient for proving computability of the Radon-Nikodym derivative $h=\frac{\rm d\mu}{\rm d\lambda}$~\cite{HRW11b}.

\section{Computability via Representations}\label{secb}

In this section we outline very shortly some concepts of the representation approach to computable analysis (TTE) \cite{Wei00,BHW08}. In TTE  computability of functions on ``concrete'' sets such as $\IN$, $\IN^\IN$, $\s$ (finite words) and $\om$ (infinite sequences) for some finite alphabet $\Sigma$ with $\{0,1\}\In\Sigma$ is defined canonically. Then these ``natural'' sets are used as sets of ``names'' of abstract objects.

A representation of a set $M$ is a surjective partial function $\delta\pf Y\to M$ where $Y$ is a natural set, $(M,\delta)$ is called a represented space.
Every $p\in\dom(\delta)$ such that $\delta(p)=x$ is called a {\em $\delta$-name} of $x$ (or {\em name} of $x$ if $\delta$ is clear from the context).
An element $x\in M$ is {\em computable} if it has a computable name.
Let $\delta_i\pf Y_i\to M_i$ ($i=1,2$) be representations.
 A realizer for a (partial) function $f\pf M_1\to M_2$ is a (partial) function $F\pf Y_1\to Y_2$ such that $f\circ\delta_1(y)=\delta_2\circ F(y)$ for all $y\in\dom(f\circ\delta_1)$.
A realizer operates on names. The function $f$ is $(\delta_1,\delta_2)$-computable if it has a computable realizer.
If $\delta_1,\delta_2$ are clear from the context we will omit the prefix and  simply say {\em computable}.  The image of a computable element by a computable function is computable.
Of course, computability on a set $M$ depends crucially on the representation.
For many spaces we use standard representations, for example $\nu_\IN$ for the natural numbers, $\nu_\IQ$ for the rational numbers and $\rho$ for the real numbers \cite{Wei00}. For computable topological spaces there are very natural standard representations, which are {\em admissible} \cite{WG09,Wei00, Sch02b}.

Consider fixed represented spaces. For every theorem of the form $(\forall x)(\exists y)Q(x,y)$ we can ask whether there is a computable (w.r.t. the given representations) function (or multi-function) mapping every $x$ to some $y$ such that $Q(x,y)$. Very often problems of this kind have no computable solution. But sometimes a solution of one problem can help to solve another one. Such ``helping'' between in general non-computable problems can be formalized by a reducibility relation as follows: For functions $f,g$ between represented spaces define  $f\leq_W g$ if there are computable (partial) functions $K,H$  on natural sets such that for every realizer $G$ of $g$, $p\mapsto K(p,G\circ H(p))$ is a realizer of $f$ ( \cite{BG11,BG11a}, where $\leq_W$ is called Weihrauch reducibility).
In other words, $f\leq_W g$ if $f$ can be computed using one single application of $g$ (provided by an oracle) in the computation   \cite[Theorem~7.2]{TW11a}. If there are  computable functions $K,H$ such that $K\circ G\circ H$ is a realizer of $f$ if $G$ is  a realizer of $g$ then $f$ is called strongly reducible to $g$, $f\leq_{sW} g$ \cite{BG11}.
These reducibilities have already been used for comparing a number of non-computable mathematical theorems, for example, in \cite{Wei92a,Wei92c,Bra05,Myl06,GM09,BG11,BG11a,BBP10}. For proving $\leq_W$, \cite[Theorem~7.2]{TW11b} is a useful tool.
The analogy of this project to reverse mathematics \cite{Sim99} has been discussed in \cite{GM09,BG11}.

An important non-computable operator is EC, the operator transforming every enumeration of a set of natural numbers into its characteristic function.
Define representations $\rm En$ and $\rm Cf$ of $2^\mathbb{N}$ by
\begin{eqnarray*}
\dom({\rm En})&=& \{0^{n_0}10^{n_1}10^{n_2}1\ldots\mid n_i\in\IN,\ (n_i=n_j\neq 0\Longrightarrow i=j)\}\\
{\rm En}(p) & = & \{n\in\IN\mid 10^{n+1}1 \mbox{ is a subword of } p\} \\
{\rm Cf}(p) & = & \{n\in\IN\mid p_n=1\}
\end{eqnarray*}
where $p\in\Sigma^\IN$ and $p_n$ is the $n$-th symbol in $p$.
Thus, if ${\rm En}(p)=A$ then $p$ enumerates $A$ without repetitions using $n_i=0$ as a dummy. The ${\rm En}$-computable sets are the r.e. sets and the $\rm Cf$-computable sets are the recursive sets.
We define the operator $\rm EC$  as the identity from the represented space $(2^\IN,{\rm En})$ to the represented space $(2^\IN,{\rm Cf})$. It transforms every enumeration of a set into its characteristic function. EC is not computable.
In \cite{Bra05} it is  proved that $\rm EC$ is equivalent to the ordinary limit map of any (sufficiently rich) computable metric space $X$. Furthermore,
EC is complete for effectively $\Sigma^0_2$-measurable functions (in the Borel hierarchy) with respect to $\leq_W$ \cite[Theorem~7.6]{Bra05}. Many non-computable problems from Analysis are equivalent to EC \cite{Wei92a,BG11,BG11a}.

\section{Computable Measurable Spaces}\label{sece}

In this section we briefly recall some basic definitions and facts from measure theory and then introduce the computational background for studying computability of the Radon-Nikodym theorem.  For a complete treatment of measure theory see for example \cite{Hal74,Bau72,Rud87,Dud02,Bog07-1}.
A~ring $\ring$ over a set $\Omega$ is a collection of subsets of $\Omega$ such that $\emptyset\in\ring$ and $ A\cup B, A\cap B, A\setminus B\in \ring$ if $A,B\in\ring$.
 A  $\sigma$-algebra $\algebra$ (over the set $\Omega$) is a collection of subsets of $\Omega$ which contains $\Omega$ and is closed with respect to  complementation and countable union. For a ring $\ring$ let $\sigma(\ring)$ be the smallest $\sigma$-algebra over $\Omega$ containing $\ring$ (the $\sigma$-algebra generated by $\ring$).
In this article we will work with a fixed measurable space $(\Omega,\algebra)$, where $\algebra$ is a $\sigma$-algebra generated by a countable ring $\ring$. Members of $\algebra=\sigma(\ring)$ will be referred to as measurable sets.

A measure on a collection $\mathcal{C}$ (which is closed  with respect to finite union) of subsets of $\Omega$ is a function $\mu:{\mathcal{C}}\to \IR^\infty$ ($=[0, \infty]$) such that i) $\mu(\emptyset) =0$, $\mu(E)\geq 0$ for all $E\in\mathcal{C}$, and ii) $\mu(\bigcup_iE_i)=\sum_i\mu(E_i)$ for pairwise disjoint sets $E_0,E_1,\ldots\in\mathcal{C}$ such that $\bigcup_iE_i\in \mathcal{C}$.
A measure $\mu$ on a collection $\mathcal{C}$ is  $\sigma$-finite, if there are sets $E_0,E_1,\ldots\in\mathcal{C}$ such that $\mu(E_i)<\infty$ for all $i$ and $\Omega=\bigcup_iE_i$.
It is well-known that every $\sigma$-finite measure on a ring $\ring$ has a unique extension to a measure on the $\sigma$-algebra $\sigma(\ring)$.

For a measure $\mu$ on a measurable space $(\Omega,{\mathcal A)}$, \, ${\mathcal L}^1(\mu)$ denotes the set of $\mu$-integrable functions $f:\Omega\to\IR$, that is, functions $f$ for which $\int f\rmd \mu$ exists. ${\mathcal L}^1(\mu)$ is an $\IR$-vector space with a seminorm $\|f\|_\mu=\int|f|\rmd\mu$. With the equivalence relation
($f\equiv_\mu g\iff \|f-g\|_\mu =0$) we obtain the normed vector space
\[L^1(\mu):={\mathcal L}^1(\mu)/\hspace{-.7ex}\equiv_\mu\hspace{1ex} =  \{[f]\mid f\in {\mathcal L}^1(\mu)\}\
\mbox{ where } \ [f]=\{g\in {\mathcal L}^1(\mu)\mid \|f-g\|_\mu=0\}\]
with the norm $\|\,[f]\,\|_\mu:=\|f\|_\mu$.
As usual in measure theory for $f:\Omega\to\IR$ we will often say
``$f\in L^1(\mu)$'' instead of ``$f\in {\mathcal L}^1(\mu)$'' or ``$[f]\in L^1(\mu)$''.

For studying computability we introduce effective versions of the definitions. In  \cite{WW06} {\em computable measure spaces} are defined. Since in this article we want to study the set of all measures on a fixed measurable space, we first define {\em computable measurable spaces} by omitting the measure in the definition from~\cite{WW06}.

\begin{definition}\label{d1}
A {\em computable measurable space} is a tuple $(\Omega,\algebra,\ring,\alpha)$ where
\begin{enumerate}
\item $(\Omega,\algebra)$ is a measurable space, $\ring$ is a countable ring such that $\Omega=\bigcup\ring$ and \linebreak $\algebra=\sigma(\ring)$,
\item $\alpha:\IN\to\ring$ is a numbering such that the operations $(A,B)\mapsto A\cup B$,
$(A,B)\mapsto A\cap B$ and $(A,B)\mapsto A\setminus B$ are computable w.r.t. $\alpha$.
\end{enumerate}
\end{definition}

\noindent {\bf In the following,  {\boldmath $(\Omega,\algebra,\ring,\alpha)$} will be a fixed computable measurable space.}
\medskip

We will consider only measures  $\mu:\algebra\to[0,\infty]$ such that $\mu(E)<\infty$ for every $E\in \ring$. Observe that on our computable measurable space such a measure $\mu$  is $\sigma$-finite, and therefore well-defined by its values on the ring $\ring$. 
We generalize the definition of a computable measure in~\cite{WW06} to representations of measures such that a measure in~\cite{WW06}  is computable if it has a computable name.

\begin{definition}\label{d2}
Let ${\mathcal M}$ be the set of measures $\mu$ such that $\mu(E)<\infty$ for all $E\in\ring$ and let ${\mathcal M}_{<\infty}$ be the set of all finite measures.
Define representations $\delm\pf\Sigma^{\IN}\to {\mathcal M}$ and
$\delmf\pf\Sigma^{\IN}\to {\mathcal M}_{<\infty}$ as follows:
\begin{enumerate}
\item
$\delm(p)=\mu$, iff $p$ is (more precisely, encodes) a list of all
$(u,n,v)\in \IQ \times\IN\times\IQ$ such that $u<\mu\circ \alpha(n)<v$,
\item
$\delmf(p)=\mu$, iff $p=\langle p_1,p_2\rangle$ such that
$\delm(p_1)=\mu$ and $p_2$ is (more precisely, encodes) a list of all
$(u,v)\in \IQ^2$ such that $u<\mu(\Omega)<v$.
\end{enumerate}
\end{definition}

The definition of $\delm$ is not artificial but follows from a very general principle. A $\delm$-name $p$ of a measure $\mu$ is a list of all (names of) $(a,E,b)\in\IQ\times \ring\times \IQ$ such that $a<\mu(E)<b$. According to \cite[Definition~8, Theorem~9]{WG09} $\delm$ is the canonical representation of the {\em effective predicate space}
$({\mathcal M}, \sigma,\nu)$, where $\nu(a,E,b)=\{\mu\in {\mathcal M}\mid a<\mu(E)<b\}$ such that $\delm$ is admissible w.r.t. the topology on ${\mathcal M}$ generated by the set $\range(\nu)$ as a subbase.

A $\delm$-name $p$ allows to compute $\mu(E)$ for every ring element $E$ with arbitrary precision. A $\delmf$-name allows additionally to compute $\mu(\Omega)$. Obviously, $\delmf\leq \delm$ and
$\delmf\equiv \delm$ if $\Omega\in\ring$.
But in general, not even the restriction of $\delta_{\mathcal M}$ to the finite measures is reducible to $\delmf$. Also, there is  a $\delm$-computable finite measure that is not $\delmf$-computable. (There is a computable sequence $(a_0,a_1,\ldots)$ of rational numbers such that $\sum_ia_i$ is finite but not computable~\cite{Wei00}. Let $\Omega:=\IN$,
${\mathcal A}:= 2^\IN$, ${\mathcal R}:=$ the set of finite subsets with canonical numbering $\alpha$, $\mu(\{n\}):= a_n$.)
  A computable measurable space with $\delta_{\mathcal M}$-computable measure $\mu$ is the ``computable measure space'' from~\cite{WW06}.
Definitions~\ref{d1} and \ref{d2} generalize the definitions from \cite{Wei99a} (probability measures on Borel subsets of $[0,1]$),
\cite{Gac05,SS06,Sch07b,GHR09a,HR09} (computable finite  measures on the Borel subsets of a computable metric space) and \cite{WW06} (computable measure on a computable $\sigma$-algebra).
\smallskip

For $E\In \algebra$, let $\chi_E:\Omega\to\IR$ be the characteristic function of $E$.
A {\em rational step function} is a finite sum
\[ s:=\sum_{k=1}^ja_k\chi_{E_k}\,,\]
where $a_k\in\IQ$ and $E_k\in\ring$. Then $s\in {\mathcal L}^1(\mu)$ and the
integral of $s$ with respect to $\mu$ is
 \[\int s\rmd\mu=\int \sum _{k=1}^ja_k\chi_{E_k}\rmd\mu= \sum_{k=1}^j a_k\mu(E_k)\,.\]
Since $\|\,.\,\|_\mu$ is a norm, $d_\mu$ with $d_\mu(f,g):=\|f-g\|_\mu$ is a metric on $L^1(\mu)$.
The set ${\rm RSF}$ of rational step functions is dense in the metric space $(L^1(\mu), d_\mu)$.

Let $\widehat\alpha$ be a canonical numbering of the set ${\rm RSF}$ of rational step functions. For fixed measure $\mu$,
$(L^1(\mu),d_\mu,{\rm RSF},\widehat\alpha)$ is an effective metric space \cite[Definition~8.1.2]{Wei00}.
The Cauchy representation $\delta_\mu\pf \om\to L^1(\mu)$ for this space is defined by:
\begin{eqnarray*}
\delta_\mu(p)=f &\iff &
\left\{
\begin{array}{l}
\mbox {$p$ \ is (more formally, encodes) a sequence}\  (s_0,s_1,\ldots) \\ \mbox{of rational step functions such that $\|s_i-f\|_\mu\leq 2^{-i}$.}
\end{array}
\right.\end{eqnarray*}

Notice that
$(\forall j>i)\|s_i-s_j\|_\mu\leq 2^{-i}\Longrightarrow (\forall i)\|s_i-f\|_\mu\leq 2^{-i} \Longrightarrow (\forall j>i)\|s_i-s_j\|_\mu\leq2\cdot 2^{-i}$. \\

For a $\delm$-computable measure $\mu$, integration $f\mapsto\int f\rmd\mu$ is
$(\delta_\mu,\rho)$-computable. We need computability also in $\mu$.
For expressing this we introduce a ``universal'' multi-representation
$\xoverline{\delta}$.
Computability w.r.t. multi-representations is defined in \cite{Wei08,TW11b}.

\begin{definition}\label{d3} Define a multi-representation $\xoverline\delta$ of the set
$\xoverline L:=\bigcup_{\mu\in{\mathcal M}}L^1(\mu)$ by
\[(\forall \mu\in{\mathcal M})\;(\forall f\in L^1(\mu)) \ \
(f\in\xoverline\delta(p)\iff f=\delta_\mu(p))\,.\]
\end{definition}

\begin{lemma}\label{l1} \hfill
\begin{enumerate}
\item
The function $H_1:{\mathcal M}\times {\rm RSF}\to \xoverline L$ where
$(H_1(\mu,s)=f :\iff s=f\in L^1(\mu))$, \\is $(\delm,\widehat\alpha,\xoverline\delta)$-computable.
\item
The function $H_2:{\mathcal M}\times {\rm RSF}\to \IR$, $H_2(\mu,s):=\int s\rmd\mu$, is $(\delm,\widehat\alpha,\rho)$-computable.
\item
The function $H_3\pf{\mathcal M}\times \xoverline L\to \IR$, where
$(H_3(\mu,f)=x :\iff f\in L^1(\mu)\wedge \int f\rmd\mu=x)$ \\
 is $(\delm,\xoverline\delta,\rho)$-computable.
\item
The function $H_4\pf{\mathcal M}\times \xoverline L\to \IR$, where
$(H_4(\mu,f)=x:\iff f\in L^1(\mu)\wedge \|f\|_\mu=x)$ \\
 is $(\delm,\xoverline\delta,\rho)$-computable.
\end{enumerate}
\end{lemma}

\pproof Straightforward.
\qq

\section{The upper bound}\label{secc}

In this section we prove that for every computable measurable space $(\Omega,\algebra,\ring,\alpha)$ the (the non-computable) operator ${\rm EC}$, mapping every enumeration of a subset of $\IN$ to its characteristic function, is an upper bound in the $\leq_W$-hierarchy of the Radon-Nikodym operator.

\begin{theorem}\label{t1} The function $\rm RN$ mapping every $\sigma$-finite measure $\lambda\in{\mathcal M}$ and every finite measure $\mu\in{\mathcal M}_{<\infty}$ such that $\mu\ll \lambda$ to the function $h\in L^1(\lambda)$ such that $\mu(E)=\int _E h\rmd\lambda$ for all $E\in{\sigma(\ring)}$, is computable via the representations $\delm$, $\delmf$ and $\xoverline\delta$ with a single application of the operator $\rm EC$. Formally, ${\rm RN}\leq_W {\rm EC}$.
\end{theorem}

\pproof Let $\lambda=\delm(p)$ and $\mu=\delmf(q)$. From $p$ and $q$ we want to find a $\delta_\lambda$-name of the function $h$.

First we ``partition'' $\lambda$ into a sequence of finite measures. Let $F_0:=\alpha(0)$ and $F_{n+1}:=\alpha(n+1)\setminus(F_0\cup\ldots\cup F_n)$. Then
$F_i\cap F_j=\emptyset$ for $i\neq j$ and $\Omega=(F_0\cup F_1\cup\ldots)$. By Definition~\ref{d1} the function $i\mapsto F_i$ is $(\nu_\IN,\alpha)$-computable and by Definition~\ref{d2}, the function $(\lambda,i)\mapsto \lambda(F_i)$ is $(\delm,\nu_\IN,\rho)$-computable. There is a $(\delm,\nu_\IN,\nu_\IN)$-computable function $d$ such that $d(\lambda,i)>\lambda(F_i)\cdot 2^i$.

For $\lambda\in{\mathcal M}$ define the function $w_\lambda:\Omega\to\IR$ by
$w_\lambda(x):=1/d(\lambda,i)$ if $x\in F_i$. Then
\[\int w_\lambda\rmd \lambda =\sum_i\lambda(F_i)/d(\lambda,i)<2\,.\]
Define a new measure $\nu$ by its values on the ring $\ring$:
\begin{eqnarray}\label{f6}\nu(E):=\int_E w_\lambda\rmd \lambda    \ \ \ (\mbox{for}\ \ E\in\ring)\,.
\end{eqnarray}
Since $w_\lambda(x)> 0$ for all $x$ and $\int w_\lambda\rmd \lambda<2$, $\nu(E)$ is well-defined for all $E\in\ring$. Since $\nu$ is a $\sigma$-finite measure on $\ring$,  it has a unique extension to the algebra $\algebra=\sigma(\ring)$. Since the function $(\lambda,E)\mapsto \nu(E) $ is $(\delm,\alpha,\rho)$-computable and the function $\lambda\mapsto \nu(\Omega)=\int w_\lambda\rmd \lambda$ is $(\lambda,\rho)$-computable, the measure $\nu$ is finite and the function $\lambda\mapsto \nu$ is $(\delm,\delmf)$-computable.

Since $w_\lambda(x)> 0$ for all $x$, $(\nu(E)=0\iff \lambda(E)=0)$ for all $E\in\ring$, hence $(\nu(A)=0\iff \lambda(A)=0)$ for all $A\in\algebra$.
Therefore, $\nu\ll\lambda\ll \nu$.

For the finite measures $\mu$ and $\nu$ with $\mu\ll\nu$ by the classical Radon-Nikodym theorem there is a function $h'\in L^1(\nu)$ such that
\begin{eqnarray}\label{f1}\mu(E)=\int_E h'\rmd\nu   \ \ \ (\mbox{for all}\ \ E\in\ring)\,.
\end{eqnarray}

First, we compute a sequence $(t'_0,t'_1,\ldots)$ of rational step functions converging to $h'$ in norm $\|\,.\,\|_\nu$.
For a partition ${\mathcal P}\In\algebra$ of $\Omega$ and a function $g\in L^1(\nu)$,
define $\mathbb E(g|{\mathcal P}):\Omega\to\IR$ as follows: for $x\in E\in{\mathcal P}$ let
\[\mathbb E(g|{\mathcal P})(x):=\left\{
\begin{array}{lll}
\int_E g\rmd\nu/\nu(E) &\mbox{if}  &\nu(E)\neq 0\\
0&         \mbox{if}  &\nu(E)= 0
\end{array}
\right.\]
 Let $ {\mathcal P_0}\In{\mathcal P}_1\In\ldots \In \algebra$ be a sequence of partitions such that  $\sigma(\bigcup_n{\mathcal P_n})=\algebra$. Then
 by Levy's zero-one law \cite[Section~10.5]{Dud02}
\begin{eqnarray}\label{f5}
\mathbb{E}(g|{\mathcal P}_n)\underset{n\to\infty}{\longrightarrow} g \ \ \ \mbox{in}\ \ L^1(\nu)\,.
\end{eqnarray}

Let $\alpha$ be the numbering of the ring $\ring$. For every $n\in\IN$ define
\[\mathcal{Q}_n:=\left\{
\begin{array}{ll}\{G_0\cap\ldots \cap G_n\mid G_j\in\{\alpha(j), \, (\alpha(0)\cup\ldots\cup\alpha(n))\setminus \alpha(j)\}\}\\
\cup  \bigcup_{i\geq n}\left(\alpha(i+1)\setminus (\alpha(0)\cup\ldots\cup\alpha(i))\right)\,.
\end{array}
\right.\]
Then $\mathcal {Q}_n\In\ring$, \, $\bigcup\mathcal{Q}_n=\bigcup_i\alpha(i)=\Omega$ and $A\cap B=\emptyset$ for $A,B\in\mathcal{Q}_n$ with $A\neq B$. Therefore, ${\mathcal P}_n:={\mathcal Q}_n\setminus\{\emptyset\}\In\ring$ is a partition of $\Omega$.
Obviously $\mathcal P_n\In \mathcal P_{n+1}$ and $\sigma(\bigcup_n{\mathcal P_n})=\algebra$.
Since union, intersection and difference on $\ring$ are computable w.r.t. $\alpha$, there is a computable function $H:\IN\times\IN\to\IN$ such that for all $n \in\IN$,
\[\mathcal{Q}_n=\{\alpha\circ H(n,i)\mid i\in\IN\} \ \ \mbox{and }\ \
\alpha\circ H(n,i)\cap \alpha\circ H(n,j)=\emptyset \ \mbox{ for }\ i\neq j\]

By (\ref{f1}), for $h_n:=\mathbb{E}(h'|{\mathcal P}_n)$ we have: for $x\in E\in {\mathcal P}_n$
\[h_n(x):=\left\{
\begin{array}{lll}
\mu(E)/\nu(E) &\mbox{if}  &\nu(E)\neq 0\\
0&         \mbox{if}  &\nu(E)= 0\,.
\end{array}
\right.\]

From ($\delmf$-names of) $\nu$ and $\mu$ and (an $\alpha$-name of) $E$ we can compute  $\mu(E)/\nu(E)$ provided $\nu(E)\neq 0$. Since $\nu(E)>0$ cannot be decided we cannot compute the functions $h_n$ ``directly''. Therefore, we
proceed as follows.
For every $n$ there are some finite set $I\In \IN$ and rational numbers $a_i$, $i\in I$, such that
\begin{eqnarray}\label{f2}\left.
\begin{array}{ll}
(\forall i\in I)\,\nu\circ \alpha\circ H(n,i)>0\,,\\
\sum _{i\in I}\mu\circ \alpha\circ H(n,i) > \mu(\Omega)-2^{-n-1}\,,\\
\sum _{i\in I}\left|\mu\circ \alpha\circ H(n,i) -a_i \cdot\nu\circ \alpha\circ H(n,i)\right| <2^{-n-1}\,.
\end{array}
\right\}\end{eqnarray}
Since the predicate (\ref{f2})  is r.e. in the variables, for $\nu,\mu$ and $n$ we can compute some finite $I\In \IN$ and rational numbers $a_i$ ($i\in I$) such that (\ref{f2}) is true. Therefore, we can compute an $\widehat \alpha$-name of $t'_n:=\sum_{i\in I}\chi_{a_i}\cdot \alpha\circ H(n,i)$. Then
by (\ref{f2}),

\begin{eqnarray*}
\|h_n-t'_n\|_\nu&=& \sum_{i\in I}\int_{\alpha\circ H(n,i)}|h_n-t'_n|\rmd\nu
+ \sum_{i\not\in I}\int_{\alpha\circ H(n,i)}|h_n-t'_n|\rmd\nu \\
&=&\sum_{i\in I}\left| \frac{\mu\circ\alpha\circ H(n,i)}{\nu\circ\alpha\circ H(n,i)}  -a_i   \right|\nu\circ\alpha\circ H(n,i)
+\sum_{i\not\in I}\int_{\alpha\circ H(n,i)}h_n\rmd\nu\\
& = & \sum_{i\in I}| \mu\circ\alpha\circ H(n,i)-a_i\cdot
\nu\circ\alpha\circ H(n,i)|
+\sum_{i\not\in I}\mu\circ \alpha\circ H(n,i) \\
&<&2^{-n}\,.\end{eqnarray*}
Since $\|h'-t'_n\|_\nu\leq \|h'-h_n\|_\nu + \|h_n-t'_n\|_\nu$, the sequence $(t'_n)_n$ of rational step functions converges to $h'$ in $L^1(\nu)$-norm. As we have observed
($\widehat\alpha$-names of) the functions $t'_n\in{\rm RSF}$ can be computed  from ($\delmf$-names of) $\mu$ and $\nu$.

Suppose $E\In F_i$ and $g\in L^1(\nu)$.
By (\ref{f6}),  $\nu(E)=\lambda(E)/d(\lambda,i)$,
hence
\begin{eqnarray}
\int_Eg\rmd\nu=\int_Eg/d(\lambda,i)\rmd\lambda=\int_Egw_\lambda\rmd\lambda\,,\\
\int g\rmd\nu=\sum_i\int_{F_i}g\rmd\nu=\sum_i \int_{F_i}gw_\lambda\rmd\lambda=\int gw_\lambda\rmd\lambda
\end{eqnarray}
Let $h:= h'w_\lambda$. Then for $E\in\ring$,
\[\mu(E)=\int_Eh'\rmd\nu=\int\chi_E h'\rmd\nu=\int\chi_E h'w_\lambda\rmd\lambda
=\int_Eh\rmd\lambda\]
Since $\|t'_nw_\lambda-h\|_\lambda=\int|t'_nw_\lambda-h'w_\lambda|\rmd\lambda  =\int|t'_n-h'|w_\lambda\rmd\lambda =\int|t'_n-h'|\rmd\nu =\|t'_n-h'\|_\nu $, the sequence $t'_nw_\lambda$ of rational step functions converges to $h$ in $L^1(\lambda)$-norm.

Let $t_n:=t'_nw_\lambda$. Then from (a $\delmf$-name of) $\mu$ and (a $\delm$-name of) $\lambda$,  ($\widehat\alpha$-names of) the $t_n\in{\rm RSF}$ can be computed such that the sequence $(t_n)_{n\in\IN}$ converges to the Radon-Nikodym derivative $h$ of $\mu$ w.r.t. $\lambda$ in $L^1(\lambda)$-norm.

We apply the non-computable operator $\rm EC$ for finding a $\xoverline\delta$-name of $h$, that is, a sequence $(s_n)_n$ of rational step functions rapidly converging to $h$ in $L^1(\lambda)$-norm.
For $\lambda\in{\mathcal M}$ and a sequence $(t_i)_i$ of rational step functions define $T\In\IN$ by
\[T:=\{\langle n,k\rangle\mid (\exists i>n)\, \|t_n-t_i\|_\lambda >2^{-k}\}\,.\]
An enumeration of $T$ can be computed from $\lambda\in{\mathcal M}$ and  $(t_i)_i$. The operator $\rm EC$ produces the characteristic function of $T$.
Let $n_0$ be the smallest number $n$ such that $\langle n,0\rangle\not\in T$
and, inductively, let $n_k$ be the smallest number $n>n_{k-1}$ such that
$\langle n,k\rangle\not\in T$. Obviously, $\| t_{n_k}-h\|_\lambda\leq 2^{-k}$.

Therefore for the computable measurable space   $(\Omega,\algebra,\ring,\alpha)$, a $\delta_\lambda$-name of the Radon-Nikodym derivative $\frac{\rmd\mu}{\rmd\lambda}$ can be computed from $\mu$ and $\lambda$ with a single application of the non-computable operator $\rm EC$.
\qq

\section{The lower bound}\label{secd}
We construct a computable measurable space $\mathcal S$ such that the operator $\rm EC$ can be computed with a single application of the Radon-Nikodym operator  for $\mathcal S$.

Let $\Omega:=[0;1)\In\IR$.
For $n\in \IN$ and $0\leq k<2^n$ let $J_{nk}:=[k/2^n;(k+1)/2^n)$ be the $k$th simple binary subinterval of length $2^{-n}$.
For $n\geq 0$ let ${\rm BI}_n:=\{J_{nk}\mid  0\leq k< 2^n\}$,
and let
${\rm BI}:=\bigcup_n{\rm BI}_n$ be the set of all simple binary intervals with canonical notation $\alpha'$. Let $\mathcal R$ be the set of finite unions of intervals from $\rm BI$ with canonical notation $\alpha$. Then ${\mathcal S}:=(\Omega, \sigma(\ring),{\mathcal R},\alpha)$ is a computable measurable space.

We will show that the operator $\rm EC$ can be computed with a single application of the Radon-Nikodym $\rm RN$ operator on $\mathcal S$. It suffices to apply the restriction of  $\rm RN$ to a fixed measure $\lambda_0$ and the finite measures $\mu\leq \lambda_0$.

Let $\lambda_0$ be the restriction of he Lebesgue measure to $\Omega:=[0;1)\In\IR$, that is $\lambda_0(J)={\rm length}(J)$ for every interval $J\In [0;1)$.
Let ${\rm RN}_0$ be the operator mapping every measure $\mu$ such that $\mu(\Omega)=1/2$ and $\mu\leq \lambda_0$ to the Radon-Nikodym derivative $h\in L^1(\lambda_0)$ such that $\mu(E)=\int_E h\rmd\lambda_0$ for all $E\in \ring$.

\begin{theorem}\label{t2}
The operator $\rm EC$ can be computed with one application of ${\rm RN}_0$, precisely,
${\rm EC}\leq_{sW} {\rm RN}_0$.
\end{theorem}

\pproof
It suffices to find computable functions $K,H$ on $\om$ such that $p\mapsto K\circ G\circ H(p)$ realizes $\rm EC$ if $G$ realizes ${\rm RN}_0$.
First, for $p\in\dom({\rm En})$ we define a function $h_p\in L^1(\mu)$.

For every interval $[a;b)\in \rm BI$ define the double jump $r[a;b)\in \rm RSF$ by
\[r[a;b)(x):=\left\{
\begin{array}{cll}
0 & \mbox{if} & x<a \mbox{ or } x\geq b\,,\\
1/2  & \mbox{if} & a\leq y<(a+b)/2\,,\\
-1/2& \mbox{if} & (a+b)/2\leq x<b\,.
\end{array}
\right.\]

Suppose $p=(0^{n_0}10^{n_1}10^{n_2}1\ldots)\in\dom({\rm En})$.
Define functions $t_k\in {\rm RSF}$ for $k\in\IN$ and $h_p:[0;1)\to\IR$ as follows:
\begin{eqnarray*}
t_k&:=& \left\{\begin{array}{ll}
0& \mbox{if}\ \ n_k=0\,,\ \\
\sum \{r(I)\mid I\in {\rm BI}_{n+1+k},\ \ I\In [2^{-n-1};2^{-n})\}& \mbox{if}\ \ n_k=n+1\,,
\end{array}
\right.\\
h_p&:=& 1/2 +\sum_{k\in\IN}t_k\,.
\end{eqnarray*}
Figure~\ref{fig1} shows the function $t_3+1/2$ for the case  $n_3=n+1$. There are $8=2^3$ up-down-up jumps in the interval $[2^{n-1};2^n)$.
\setlength{\unitlength}{0,72pt}
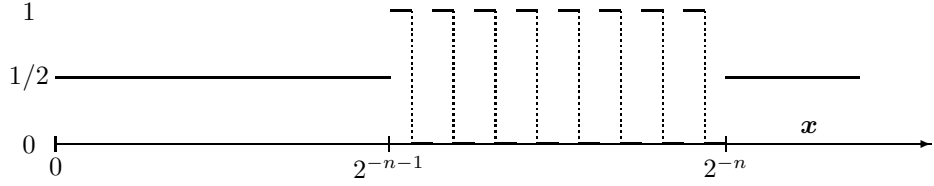
\begin{figure}[htbp]

\newsavebox{\doublejump}
\thicklines
\savebox{\doublejump}(0,0)[]{
\put(0,70){\line(1,0){11}}
\put(11,0){\line(1,0){11}}
\thinlines
\multiput(11,0)(0,3){24}{\line(0,1){1.0}}
}
\begin{picture}(320,100)(-90,-20)
\thinlines
\put(-76,0){\vector(1,0){460}}
\put(-76,4){\thicklines \line (0,-1){8}}
\put(99,4){\line (0,-1){8}}
\put(276,4){\line (0,-1){8}}
\put(320,9){\makebox(0,0){\boldmath $x$}}

\put(-76,-12){\makebox(0,0){$0$}}
\put(99,-12){\makebox(0,0){$2^{-n-1}$}}
\put(276,-12){\makebox(0,0){$2^{-n}$}}
\put(-90,70){\makebox(0,0){$1$}}
\put(-90,35){\makebox(0,0){$1/2$}}
\put(-90,0){\makebox(0,0){$0$}}

\thicklines
\put(-76,35){\line (1,0){176}}
\put(276,35){\line (1,0){70}}

\multiput(100,35)(22,0){8}{\usebox{\doublejump}}
\thinlines

\end{picture}
\caption{The function $t_3+1/2$ for the case  $n_3=n+1$.} \label{fig1}
\end{figure}
\\
If $A={\rm En}(p)$ then $h_p$ has up-down-up jumps in the interval $[2^{-n-1};2^{-n})$ iff $n\in A$, and the number of jumps in this interval is $2^k$ where $n_k=n+1$.

First, we observe that for all $I,I'\in{\rm BI}$, $I\cap I'=\emptyset$, $ I\In I'$ or $I'\In I$. Therefore
\[\int_E r(I)\rmd\lambda_0=0 \ \ \mbox{for all}\ \
E,I\in{\rm BI} \ \ \mbox{such that} \ \  {\rm length}(I)\leq {\rm length}(E))\,.
\]
Then for $E\in {\rm BI}_m$,
\begin{eqnarray*}
&&\int_E\sum_{k\in\IN}t_k\rmd\lambda_0 \\
&=& \int_E\sum _{k\in\IN}
\sum \{r(I)\mid I\in {\rm BI}_{n+1+k},\ \ I\In [2^{-n-1};2^{-n}),\ n_k=n+1\}\rmd\lambda_0\\
&=& \sum _{k\in\IN}
\sum \{\int_E r(I)\rmd\lambda_0\mid I\in {\rm BI}_{n+1+k},\ \ I\In [2^{-n-1};2^{-n}),\ n_k=n+1\}\\
&=& \sum _{k<m}
\sum \{\int_E r(I)\rmd\lambda_0\mid I\in {\rm BI}_{n+1+k},\ \ I\In [2^{-n-1};2^{-n}),\ n_k=n+1\}
\end{eqnarray*}
which is a finite sum.
Therefore, from $p$ and (any $\alpha$-name of) $E$ we can compute $\int_Eh_p\rmd\lambda_0$. This is true also for all $E\in\ring$.
Since $h_p\geq 0$ and $h_p\in L^1(\lambda_0)$ a measure $\mu_p$ on the computable measurable space $\mathcal S$ is defined by $\mu_p(E):=\int_E h_p\rmd\lambda_0$.
Therefore, there is a computable function $H\pf\om\to\om$ that maps every $p\in\dom({\rm En})$ to some $\delmf$-name of $\mu_p$.

Since $h_p\leq 1$, $\mu_p\leq\lambda_0$. Obviously, $h_p$ is the Radon-Nikodym derivative of $\mu_p$ w.r.t. $\lambda_0$ and $\mu_p(\Omega)=1/2$. Let $G$ be an arbitrary $(\delmf,\delta)$-realizer of the operator ${\rm RN}_0$. Then for   every
$p\in \dom({\rm En})$,
$G\circ H(p)$ is a $\delta_{\lambda_0}$-name of $h_p$.
From the definition of $h_p$,
\begin{eqnarray}\label{f7}\int_{[2^{-n-1};2^{-n})} |h_p-1/2|\rmd\lambda_0=\left\{
\begin{array}{ll}
2^{-n-1}/2 &\mbox{if} \ \  n\in {\rm En}(p)\\
0&\mbox{otherwise}
\end{array}
\right.
\end{eqnarray}
(see Figure~\ref{fig1}). The operator $f\mapsto |f-1/2|$ is $(\delta_{\lambda_0},\delta_{\lambda_0})$-computable, hence by Lemma~\ref{l1}, the function $(f,E)\mapsto\int_E|f-1/2|\rmd\lambda_0$ is
$(\delta_{\lambda_0},\alpha,\rho)$-computable.
Therefore by (\ref{f7}) there is a computable function $K\pf \om\to\om$  that from a $\delta_{\lambda_0}$-name of $h_p$ computes a ${\rm Cf}$-name of ${\rm En}(p)$.

In summary, there are computable functions $H,K$ on $\om$ such that
$K\circ G\circ H$ realizes $\rm EC$ if $G$ realizes $\rm RN_0$, hence   ${\rm EC}\leq_{sW}{\rm RN}_0$.
\qq

\section{Final remarks}\label{secf}
In our proof of the effective version of the Radon-Nikodym theorem, Theorem~\ref{t1}, via the function $w_\lambda$ we can replace the (possibly) infinite measure $\lambda$ by the finite measure $\nu$. For solving the problem for finite measures $\mu\ll\nu$, we have assumed that the classical Radon-Nikodym theorem is true and have applied
 Levy's zero-one law \cite[Section~10.5]{Dud02} for showing that our sequence $(h_n)_n$ converges to the Radon-Nikodym derivative $h'=\frac{\rmd\mu}{\rmd\nu}$. Elstrod \cite[Page 278]{Els99} gives a simple elementary proof for finite measures $\mu\leq \nu$ based on \cite{Bra89}, which essentially includes the proof of   Levy' zero-one law for this case and the definition of a sequence of step functions converging to $\frac{\rmd\mu}{\rmd\nu}$.
This proof can be easily turned into an effective version, where at the end, the operator $\rm EC$ must be applied to produce a fast converging sequence. Then the computability result can be extended from finite measures $\mu\leq  \nu$ to finite measures $\mu\ll\nu$ and to finite $\mu$ and $\sigma$-finite $\lambda$ such that $\mu\ll\lambda$ \cite{Bra89}.

Another proof of the classical theorem by J. v. Neumann applies the Fr\'echet-Riesz representation theorem for continuous linear functionals on Hilbert spaces using the fact that $L^2(\mu)$ with the product $(f;g):=\int fg\rmd\mu$ is a Hilbert space.
Also this proof can be effectivized. A computable version of the Fr\'echet-Riesz representation theorem for computable Hilbert spaces has been proved in \cite{BY06}:
from a linear functional $F$ and its norm a point $a_F$ can be computed such that
$F(x)=(x;a_F)$. In our application, $a_F$ will be the Radon-Nikodym derivative.
But for our uniform theorem we need a version
of the Fr\'echet-Riesz theorem
uniform on the class $\mathcal{ EH}$ of all  ``effective'' Hilbert spaces with $\rm RSF$ as a dense subspace. We can prove:
From a space $H\in\mathcal {EH}$ and the linear functional $F$ a sequence of rational step functions can be computed converging (not necessarily fast) to $a_F$. By an application of $\rm EC$ a fast converging subsequence can be selected.
Also, form $H$, $F$ and its norm $\|F\|$ a sequence of rational step functions can be computed fast converging  to $a_F$. Notice that $\|F\|$ can be computed from below (a $\rho_<$-name \cite{Wei00}) and that (a $\rho$-name \cite{Wei00} of)
$\|F\|$ can be computed from $F$ with  a single application of $\rm EC$.
In \cite{HRW11b} it is shown that for computable measures a specific computability condition (``$\mu$ is computably normable relative to $\lambda$'') on the measures $\lambda$ and $\mu$ suffices to compute a $\delta_\lambda$-name of the Radon-Nikodym derivative without using  the operator~$\rm EC$. The proof includes a direct proof of the computable
Fr\'echet-Riesz representation theorem.

\bibliographystyle{alpha}

\end{document}